\documentclass{article}
\usepackage{stmaryrd,amssymb,amsmath,amsthm,url}

\title{Differential Meadows}

\author{
	Jan A.\ Bergstra\thanks{Partially supported by the 
	  Dutch NWO Jacquard project \emph{Symbiosis}, 
	  project number 638.003.611} 
	\and
	Alban Ponse \\
\\
  {\small
	  Section Software Engineering,
	  Informatics Institute,
	  University of Amsterdam}\\
{\small URL: \url{www.science.uva.nl/~{janb,alban}}
	}
}
\date{}
\newcommand{\NM}{\ensuremath{\mathbb M}}
\newcommand{\NC}{\ensuremath{{\mathbb C}_0}}
\newcommand{\DE}{\ensuremath{\textit{DE}}}

\newcommand{\Mod}{\axname{Mod}}
\newcommand{\CRU}{\itname{CRU}}

\newtheorem{theorem}{Theorem}
  
\newtheorem{proposition}{Proposition}  
  
\newtheorem{corollary}{Corollary}  
  
\newtheorem{definition}{Definition}  
\theoremstyle{definition}

\newcommand{\IL}{\itname{IL}}

\newcommand{\ZTF}{\itname{ZTF}}

\newcommand{\Md}{\itname{Md}}
\newcommand{\Refl}{\itname{Refl}}
\newcommand{\Ril}{\itname{RIL}}

\newcommand{\itname}[1]{\ensuremath{\textit{#1}\,}}
\newcommand{\axname}[1]{\ensuremath{\textrm{#1}}}

\begin{document}

\maketitle

\begin{abstract}
A meadow is a zero totalised field ($0^{-1}=0$), and a
cancellation meadow is a meadow without proper zero divisors.
In this paper we consider differential meadows, i.e., 
meadows equipped with differentiation operators. 
We give an equational axiomatization of these
operators and thus obtain a finite basis for
differential cancellation meadows. 
Using the Zariski topology we prove the existence of a 
differential cancellation meadow.
\end{abstract}

\section{Introduction}
\label{sec:1}
A \emph{meadow} is an algebra in the signature of fields 
with an inverse operator that satisfies the equations of 
commutative rings with unit ($\CRU$) together with
\begin{align*}
(x^{-1})^{-1}&=x&\axname{(\Refl)}\\
x\cdot x\cdot x^{-1}&=x&(\Ril)
\end{align*}
where the names of the equations abbreviate \emph{Reflection} 
and \emph{Restricted Inverse Law}, respectively.
Meadows were introduced in \cite{BT07}.

In~\cite{BHT08} it was shown that the variety of meadows 
satisfies precisely those equations which are valid in all 
so-called zero totalised fields (\ZTF s). 
A \ZTF\ is a field equipped with an inverse operator
$(\_)^{-1}$ that has been made total by putting $0^{-1}=0$. 
Alternatively and following \cite{BPZ07}, we will 
qualify a zero totalised field as a \emph{cancellation meadow} 
if it enjoys the following cancellation property:
\begin{align}
\label{canc}
x \neq 0 \wedge x \cdot y &= x \cdot z 
\quad\Longrightarrow\quad
y = z.
\end{align}

The mentioned result from~\cite{BHT08} may be viewed as a 
completeness theorem: $\CRU+\Refl+\Ril$
completely axiomatises the equational theory $E(\ZTF)$ of the 
class \ZTF\ of zero totalised fields.
Another way of looking at this result is that it establishes that
$E(\ZTF)$ has a finite basis.
The proof of the finite basis theorem for $E(\ZTF)$ 
in~\cite{BHT08}
makes use of the existence of maximal ideals. Although concise and
readable,
that proof is non-elementary because the existence of maximal ideals
requires a non-elementary set theoretic principle, independent 
of ZF set theory.
In \cite{BT07} a finite basis theorem was established for $E_C(\ZTF)$, 
the closed equations true in \ZTF s. 

In~\cite{BP08} a proof of the finite basis
result for $E(\ZTF)$ has been given along the lines of the elementary
proof about $E_C(\ZTF)$.
The proof method is more general than the proof using maximal ideals
because it generalizes to extended signatures 
(see Theorem~\ref{st} below). 
In this paper we apply this result to so-called 
differential meadows, i.e., meadows equipped with formal variables
$X_1,...,X_n$ and differential operators
\[\frac{\partial}{\partial X_i}.\]
We provide a short equational axiomatization of the differential
operators and thus obtain a finite basis for
differential cancellation meadows. 
This appears to be an elegant axiomatization, e.g., 
$\frac{\partial}{\partial X_i}(1/x)=-(1/x^2)\cdot
\frac{\partial}{\partial X_i}(x)$ follows easily.
Finally, we prove the existence of a differential cancellation
meadow, using the Zariski topology~\cite{Zar44,Har77} and a representation result 
from~\cite{BP08}.

The paper is structured as follows: in the next section we 
recall cancellation meadows and the generic basis theorem, 
and introduce differential meadows.
Then, in Section~\ref{sec:3} we prove the existence of a 
differential cancellation meadow. Some conclusions 
are drawn in Section~\ref{sec:conc}. 

\section{Cancellation and Differential Meadows}
\label{sec:2}
In this section we fix some notation and explain
cancellation meadows and our generic basis theorem in detail. 
Then we introduce differential meadows.

\subsection{Cancellation meadows and a generic basis result}
A meadow is an algebra in the signature of fields that satisfies
the axioms in Table~\ref{Md}. We write \Md\ for the set of axioms 
in Table~\ref{Md}, thus (referring to the Introduction)
$\Md = \CRU+\Refl+\Ril$.

Let \IL\ (Inverse Law) stand for 
\[x\neq 0 \quad\Longrightarrow\quad x\cdot x^{-1}=1,\]
so \IL\ states that there are no zero divisors. 
Note that \IL\ and the cancellation property~\eqref{canc} 
are equivalent. 
A \emph{cancellation meadow} is a meadow that
also satisfies \IL.

\begin{table}
\centering
\hrule
\begin{align*}
	(x+y)+z &= x + (y + z)\\
	x+y     &= y+x\\
	x+0     &= x\\
	x+(-x)  &= 0\\
	(x \cdot y) \cdot  z &= x \cdot  (y \cdot  z)\\
	x \cdot  y &= y \cdot  x\\
	1\cdot x &= x \\
	x \cdot  (y + z) &= x \cdot  y + x \cdot  z\\
	(x^{-1})^{-1} &= x \\
	x \cdot (x \cdot x^{-1}) &= x
\end{align*}
\hrule
\caption{The set \Md\ of axioms for meadows}
\label{Md}
\end{table}

From the axioms in \Md\ the following identities are derivable:
\begin{align*}
	(0)^{-1}  &= 0,&0\cdot x  &= 0,\\
	(-x)^{-1} &= -(x^{-1}),&x\cdot -y &= -(x\cdot y),\\
	(x \cdot  y)^{-1} &= x^{-1} \cdot  y^{-1},&-(-x)     &= x.
\end{align*}

We write $\Sigma_m=(0,1,+,\cdot,-,^{-1})$ for the
signature of (cancellation) meadows. Furthermore,
we often write $1/t$ or 
$\displaystyle\frac 1 t $
for $t^{-1}$, $tu$ for $t\cdot u$, $t/u$ for $t\cdot (1/u)$, 
$t-u$ for $t+(-u)$, 
and freely use numerals and exponentiation with constant
integer exponents. 
We further use the notation
\[1_x\text{ for } \frac x x\quad \text{and}
\quad 0_x\text{ for }
1-1_x.\]
Note that for all terms $t$,
$(1_t)^2=1_t$, $1_t\cdot 0_t = 0$ and $(0_t)^2=0_t$. 
We call an expression $1_t$ a \emph{pseudo unit} because 
it is almost equivalent to the unit 1, and for a similar
reason we say that $0_t$ is a \emph{pseudo zero}.

The basis result from~\cite{BP08} admits generalization
if pseudo units and pseudo zeros propagate in the context rule
for equational logic. 
We recall the precise definition of this form of
propagation from that paper.

\begin{definition}
Let $\Sigma$ be an extension of $\Sigma_m=(0,1,+,\cdot,-,^{-1})$, 
the signature of meadows, and let $E\supseteq \Md$ be a set of 
equations over $\Sigma$. Then
\begin{enumerate}
\item
$(\Sigma,E)$ has the \textbf{propagation property for pseudo units} 
if for each pair of $\Sigma$-terms $t,r$ and context $C[~]$,
\[E\vdash 1_t\cdot C[r]=1_t\cdot C[1_t\cdot r].\]
\item
$(\Sigma,E)$ has the \textbf{propagation property for pseudo zeros} 
if for each pair of $\Sigma$-terms $t,r$ and context $C[~]$,
\[E\vdash 0_t\cdot C[r]=0_t\cdot C[0_t\cdot r].\]
\end{enumerate}
\end{definition} 

\noindent
We now recall our generic basis result from~\cite{BP08}:

\begin{theorem}
[Generic basis theorem for cancellation meadows]
\label{st}
If $\Sigma\supseteq \Sigma_m,~ E\supseteq \Md$ is 
a set of equations over
$\Sigma$, and $(\Sigma,E)$ 
has the
pseudo unit propagation property and the pseudo zero propagation 
property, then $E$ is a basis (a complete axiomatisation) of 
$\Mod_\Sigma(E\cup\IL)$.
\end{theorem}

\subsection{Differential Meadows}
Given some $n\geq 1$ we extend the signature $\Sigma_m$ of meadows 
with differentiation operators and constants $X_1,...,X_n$ to model
functions to be differentiated:
\[\frac{\partial}{\partial X_i}: \NM\rightarrow \NM\]
for $i=1,...,n$ and some meadow \NM. 
We write $\Sigma_{md}$ for this extended signature.
Equational axioms for $\frac{\partial}{\partial X_i}$
are given in Table~\ref{t1}, where D4 and D5 define $n^2$ 
equational axioms. 
Observe that the \Md\ axioms together with \axname{D3} imply
$\frac{\partial}{\partial X_i}(0)=0$. 
Furthermore, using axiom \axname{D1} one easily proves: 
$\frac{\partial}{\partial X_i}(-x) =-\frac{\partial}{\partial X_i}(x)$.

\begin{table}
\hrule
\begin{align*}
\frac{\partial}{\partial X_i}(x+y)&=\frac{\partial}{\partial X_i}(x)
+\frac{\partial}{\partial X_i}(y)&\axname{(D1)}\\
\frac{\partial}{\partial X_i}(x\cdot y)&=
\frac{\partial}{\partial X_i}
(x)\cdot y+x\cdot \frac{\partial}{\partial X_i}(y)&\axname{(D2)}\\
\frac{\partial}{\partial X_i}(x\cdot x^{-1})&=0&\axname{(D3)}\\
\frac{\partial}{\partial X_i}(X_i)&=1&\axname{(D4)}\\
\frac{\partial}{\partial X_i}(X_j)&=0 \quad
\text{ if } i \neq j\,&\axname{(D5)}
\end{align*}
\hrule
\caption{The set of axioms \DE}
\label{t1}
\end{table}

First we establish the expected corollary of Theorem~\ref{st}:
\begin{corollary}
\label{cor:diff}
The set of axioms $\Md\cup\DE$ (see Tables~\ref{Md} and \ref{t1})
is a finite basis (a complete axiomatisation) of 
$\Mod_{\Sigma_{md}}(\Md\cup\DE\cup\IL)$.
\end{corollary}
\begin{proof}
The pseudo unit propagation property requires a check for 
$\frac{\partial}{\partial X_i}(\_)$ only:
\begin{equation}
\label{tja}
\frac{\partial}{\partial X_i}(1_t\cdot r)=
\frac{\partial}{\partial X_i}(1_t)\cdot r+1_t\cdot 
\frac{\partial}{\partial X_i}(r)=
1_t\cdot\frac{\partial}{\partial X_i}(r).
\end{equation}
Multiplication with $1_t$ now yields the property.
From \eqref{tja} we get
\[0_t\cdot \frac{\partial}{\partial X_i}(r)
=\frac{\partial}{\partial X_i}(r)-1_t\cdot 
\frac{\partial}{\partial X_i}(r)
\stackrel{\eqref{tja}}=\frac{\partial}{\partial X_i}(r)-
\frac{\partial}{\partial X_i}(1_t\cdot r)
=
\frac{\partial}{\partial X_i}(0_t\cdot r)\]
and multiplication with $0_t$
then yields the pseudo zero propagation property.
\end{proof}

A \emph{differential meadow} is a meadow
equipped with formal variables $X_1,...,X_n$ and differentiation
operators $\frac{\partial}{\partial X_i}(\_)$ that satisfies 
the axioms in \DE.

We conclude his section with an elegant consequence of the 
fact that we are working in the setting of meadows, namely 
the consequence that the differential of an inverse follows 
from the \DE\ axioms.

\begin{proposition}
$
\displaystyle
\Md\cup\DE\,\vdash \frac{\partial}{\partial X_i}(1/x)=-(1/x^2) 
\cdot 
\frac{\partial}{\partial X_i}(x)$.
\end{proposition}

\begin{proof}
By axioms D3 and D2,\quad
$\displaystyle 0=\frac{\partial}{\partial X_i}(x/x)
=
\frac{\partial}{\partial X_i} (x)\cdot 1/x
+x\cdot\frac{\partial}{\partial X_i}(1/x)$,\quad
so 
\begin{align*}
0&=0\cdot(1/x)= \displaystyle
\frac{\partial}{\partial X_i}(x/x)\cdot(1/x)
=\frac{\partial}{\partial X_i}(x)\cdot 1/x^2+(x/x)\cdot
\frac{\partial}{\partial X_i}(1/x)\\
& \displaystyle\stackrel{\eqref{tja}}=1/x^2 \cdot 
\frac{\partial}{\partial X_i}(x)
+
\frac{\partial}{\partial X_i}((x/x)\cdot(1/x))
\stackrel{\Ril}=1/x^2 \cdot \frac{\partial}{\partial X_i}(x)+
\frac{\partial}{\partial X_i}(1/x),
\end{align*}
and hence
\[
\frac{\partial}{\partial X_i}(1/x)=-(1/x^2) \cdot 
\frac{\partial}{\partial X_i}(x).
\]
\end{proof}

\section{Existence of Differential Meadows}
\label{sec:3}
In this section we show the \emph{existence} of differential 
meadows with formal variables $X_1,...,X_n$ for arbitrary 
finite $n>0$. 
First we define a particular cancellation meadow, and then we expand 
this meadow to a differential cancellation meadow by adding formal
differentiation.

\subsection{The Zariski topology congruence over $\NC^n$}
We will use some terminoloy from algebraic geometry, 
in particular we will use the Zariski topology~\cite{Zar44,Har77}. 
Open (closed) sets in this topology will be indicated 
as Z-open (Z-closed). Recall that 
complements of Z-closed sets are Z-open and complements of 
Z-open sets are Z-closed,
finite unions of Z-closed sets are Z-closed, and
intersections of Z-closed sets are Z-closed.
Let $\NC$ denote the zero-totalized expansion of the complex numbers. 
We will make use of the following facts:

\begin{enumerate} 
\item \label{een}
The solutions of a set of polynomial equations (with $n$ or less 
variables) within $\NC^n$ constitute a Z-closed subset of $\NC^n$. 
Here 'polynomial' has the conventional meaning, not involving division. 
Taking equations $1=0$ and $0=0$ respectively, it follows that both 
$\emptyset$ and $\NC^n$ are Z-closed (and Z-open as well).

\item \label{twee}
Intersections of non-empty Z-open sets are non-empty.
\end{enumerate}
In the following we consider terms 
\[t(\overline X)=t(X_1,...,X_n)\]
with $t=t(\overline x)$ a $\Sigma_m$-term and we write 
$T(\Sigma_m(\overline X))$ for the set of these terms.
For $V\subseteq\NC^n$ we define the equivalence
\[\equiv^V_{\NC^n}\]
on $T(\Sigma_m(\overline X))$
by $t(\overline X)\equiv^V_{\NC^n}r(\overline X)$ 
if each assignment $\overline X\mapsto V$
evaluates both sides to  equal values in \NC. 
It follows immediately that for each $V\subseteq\NC^n$, 
$T(\Sigma_m(\overline X))/\equiv^V_{\NC^n}$ is a 
meadow. In particular, if 
$ V = \emptyset$ one obtains the trivial meadow ($0=1$) as both 0 
and 1 satisfy any universal quantification over an empty set. 
If $V$ is{ a singleton this quotient is a cancellation meadow. 
In other cases the meadow may not satisfy 
the cancellation property. 
Indeed, suppose that $n=1$ and $V = \{0,1\}$ and
let $t(X) = X$. Now $t(1) \neq 0$. Thus
$t(X) \neq 0$ in  $T(\Sigma_m(X))/\equiv^V_{\NC}$. 
If that is assumed to be a cancellation meadow, however, 
one has $1_{t(X)}=1$, but $1_{t(0)}=0$, thus refuting 
$1_{t(X)}=1$.

We now define the relation $\equiv_{ZTC}$ 
(Zariski Topology Congruence over $\NC^n$) by
\[t\equiv_{ZTC} r \iff \exists V (
  \text{$V$ is Z-open, $V\neq\emptyset$ and 
        $t\equiv^V_{\NC^n}r$)}.
\]
The relation $\equiv_{ZTC}$ is indeed a congruence 
for all meadow operators: 
the equivalence properties follow easily; 
for $0\equiv_{ZTC} 0$ and $1\equiv_{ZTC} 1$,
take $V=\NC^n$, and if 
$P\equiv_{ZTC} P'$ and $Q\equiv_{ZTC} Q'$, witnessed 
respectively by $V$ and $V'$, 
then $P+P'\equiv_{ZTC} Q+Q'$ and $P\cdot P'\equiv_{ZTC} Q\cdot Q'$ 
are witnessed by
$V\cap V'$  which is Z-open and non-empty because  of
fact~\ref{twee} above. Finally $-P\equiv_{ZTC} -P'$ and 
$(P)^{-1}\equiv_{ZTC}(P')^{-1}$
are both witnessed by $V$.

In ~\cite{BP08} we defined the \emph{Standard
Meadow Form} (SMF) representation result for meadow terms.
This result implies for 
$T(\Sigma_m(\overline X))/\equiv_{ZTC}$ that each term can be 
represented by 0 or by $p/q$ with $p$ and $q$ polynomials not 
equal to 0. We notice that it is decidable whether or not a 
polynomial equals the 0-polynomial by taking all corresponding 
products of powers of the $X_1,...,X_n$ together and then 
checking that all coefficients vanish. 

A few more words on the SMF representation result. SMFs are defined 
using levels: an SMF of level 0 is of the form
$p/q$ with $p$ and $q$ polynomials,
and an SMF of level $k+1$ is of the form $0_p\cdot P+1_p\cdot Q$ with
$P$ and $Q$ both SMFs of level $k$. 
As an example, let $P$ be the SMF of level 1 defined by
\[P=0_{1-X_1}\cdot \frac{2X_1}{X_2}+1_{1-X_1}\cdot 
\frac{1+X_2-2X_1X_3}{8-X_1X_3^2}.\]
Now in $T(\Sigma_m(\overline X))/\equiv_{ZTC}$, the polynomial 
$1-X_1$ is on some Z-open non-empty set $V$ not equal to 0 
(see fact~\ref{een} above), thus
$1_{1-X_1}\equiv^V_{\NC^n}1$ and $0_{1-X_1} \equiv^V_{\NC^n} 0$, 
and hence 
\[P\equiv_{ZTC}\frac{1+X_2-2X_1X_3}{8-X_1X_3^2}.\]
So, in $T(\Sigma_m(\overline X))/\equiv_{ZTC}$, 
the SMF level-hierarchy collapses and
terms can be represented by either $0$ or by $p/q$ with both $p$ and 
$q$ polynomials not equal to 0. In the second case
$1_{p/q} = 1$ and therefore it is a cancellation meadow.
Furthermore, equality is decidable in this model. 
Indeed to check that $1_p = 1$ (and $0_p = 0$) for a polynomial 
$p$ it suffices to check that $p$ is not 0 over the complex numbers. 
Using the SMF representation all closed terms are either 0 or 
take the form $p/q$ with $p$ and $q$ nonzero polynomials. 
For $q$ and $q'$ nonzero 
polynomials we find that $p/q \equiv_{ZTC} p'/q' 
\iff p\cdot q' - p'\cdot q = 0$ 
which we have already found to be decidable.

\subsection{Constructing a differential cancellation meadow}
In $T(\Sigma_m(\overline X))/\equiv_{ZTC}$ the differential 
operators can be defined as follows:
\[\frac{\partial}{\partial X_i}(0)=0\]
and, using the fact that differentials on polynomials are known,
\[\frac{\partial}{\partial X_i}(\frac p q) =
\frac{\frac{\partial}{\partial X_i}(p) \cdot 
q - p \cdot \frac{\partial}{\partial X_i}(q)}{q^2}.\]
Let $V$ be the set of 0-points of $q$ and let $U={\sim} V$, 
the complement of $V$. 
Then $p/q$ is  differentiable on $U$ and the derivative coincides 
with the formal derivative used in the definition. 
This definition is representation independent:
consider $p'/q'\equiv_{ZTC}p/q$ with $V'$ the 0-points of $q'$
and $U'={\sim}V'$. Then there is some 
non-empy and Z-open $W$
such that $p/q\equiv^W_{\NC^n}p'/q'$. Now
$W\cap U\cap U'$ is non-empty and Z-open, and on this set,
\[\frac{\partial}{\partial X_i}(\frac p q)=
\frac{\partial}{\partial X_i}(\frac {p'}{q'}).\]
So, formal differentation $\partial/\partial X_i$ preserves 
the congruence properties. Finally, we check the soundness of the
\DE\ axioms:
\begin{description}
\item[\textnormal{Axiom D1:}] 
Consider $t+t'$. In the case that one of $t$ 
and $t'$ equals 0, axiom D1 is obviously sound. 
In the remaining case,
$t=p/q$ and $t'=p'/q'$ with all polynomials
not equal to 0 and $t+t'=\frac{pq'+p'q}{qq'}$. 
Using ordinary differentiation on polynomials we derive
\begin{align*}
\frac{\partial}{\partial X_i}(t+t')
&= \frac{\frac{\partial}{\partial X_i}(pq'+p'q) \cdot 
qq' - (pq'+p'q) \cdot \frac{\partial}{\partial X_i}(qq')}{(qq')^2}\\
&=\frac{\frac{\partial}{\partial X_i}(p)\cdot q\cdot (q')^2+
\frac{\partial}{\partial X_i}(p')\cdot q^2\cdot q' 
- p\cdot \frac{\partial}{\partial X_i}(q)\cdot (q')^2
- p'\cdot \frac{\partial}{\partial X_i}(q')\cdot q^2}{(qq')^2}\\
&=\frac{\partial}{\partial X_i}(\frac{p}{q})\cdot 1_{(q')^2}+
\frac{\partial}{\partial X_i}(\frac{p'}{q'})\cdot 1_{q^2}\\
&=\frac{\partial}{\partial X_i}(t)+\frac{\partial}{\partial X_i}(t').
\end{align*}

\item[\textnormal{Axiom D2:}] Similar.

\item[\textnormal{Axiom D3:}]
Consider $t$, then either $t=0$ or $t/t=1$, and in both cases
$\displaystyle\frac{\partial}{\partial X_i}(\frac t t)=0$.

\item[\textnormal{Axioms schemes D4 and D5:}] We derive
\[\frac{\partial}{\partial X_i}(X_j)=
\frac{\partial}{\partial X_i}(\frac{X_j}1)=
\begin{cases}
0&\text{if $i\neq j$,}\\
1&\text{otherwise.}
\end{cases}\]
\end{description}

Thus, by adding formal differentiation to $T(\Sigma_m(\overline X))$
we constructed a differential cancellation meadow.

\section{Conclusions}
\label{sec:conc}
In this paper we introduced differential meadows. We provided
a finite equational basis for differential cancellation meadows and 
proved their existence by a construction based on the Zariski 
topology.

Differential meadows generalize differential fields in the same way 
as meadows generalize fields. As stated in \cite{BHT08}, exactly 
the von Neumann regular rings admit expansion to a meadow. 
The general question, however, which meadows can be expanded to 
differential meadows that satisfy the \DE\ axioms is left open. 
In \cite{BR08} finite meadows have been characterized as direct 
sums of finite fields. The existence of differential meadows 
over a finite meadow is in particular left for further analysis.

\end{document}